\documentclass[11pt]{article}
\usepackage{amscd}
\usepackage{amsmath}
\usepackage{amsfonts}
\newcommand{\al}{\alpha}

\newcommand{\dbar}{\overline{\partial}}
\newcommand{\dd}[1]{\partial_{#1}}
\newcommand{\dbr}[1]{\partial_{\overline{#1}}}

\newcommand{\ddt}[1]{\frac{\partial #1}{\partial t}}
\newcommand{\R}[4]{R_{#1 \overline{#2} #3 \overline{#4}}}

\newcommand{\lt}{\tilde{\triangle}}
\newcommand{\chii}[2]{\chi^{#1 \overline{#2}}}
\newcommand{\ch}[2]{\chi_{#1 \overline{#2}}}
\newcommand{\chz}[2]{\chi_{0 \, #1 \overline{#2}}}
\newcommand{\chp}[2]{\chi'_{#1 \overline{#2}}}
\newcommand{\g}[2]{g_{#1 \overline{#2}}}
\newcommand{\gi}[2]{g^{#1 \overline{#2}}}
\newcommand{\h}[2]{h^{#1 \overline{#2}}}
\newcommand{\tphi}{\hat{\phi}}

\begin{document}
\newcounter{remark}
\newcounter{theor}
\setcounter{remark}{0}
\setcounter{theor}{1}
\newtheorem{claim}{Claim}
\newtheorem{theorem}{Theorem}[section]
\newtheorem{proposition}{Proposition}[section]
\newtheorem{lemma}{Lemma}[section]
\newtheorem{defn}{Definition}[theor]
\newtheorem{corollary}{Corollary}[section]
\newenvironment{proof}[1][Proof]{\begin{trivlist}
\item[\hskip \labelsep {\bfseries #1}]}{\end{trivlist}}
\newenvironment{remark}[1][Remark]{\addtocounter{remark}{1} \begin{trivlist}
\item[\hskip
\labelsep {\bfseries #1  \thesection.\theremark}]}{\end{trivlist}}

\centerline{\bf ON THE CONVERGENCE AND SINGULARITIES}
\centerline{\bf OF THE J-FLOW WITH APPLICATIONS}
\centerline{\bf TO THE MABUCHI ENERGY}

\bigskip
\bigskip
\begin{tabular}{ccc}
Jian Song & & Ben Weinkove \\
Department of Mathematics & & Department of Mathematics \\
Johns Hopkins University & & Harvard University \\
3400 North Charles Street & & 1 Oxford Street \\
Baltimore MD 21218 & & Cambridge MA 02138 \\
jsong@math.jhu.edu & & weinkove@math.harvard.edu
\end{tabular}
\bigskip

\setlength\arraycolsep{2pt}
\addtocounter{section}{1}

\bigskip
\noindent
{\bf 1. Introduction}
\bigskip

The J-flow is a parabolic flow on K\"ahler manifolds with two K\"ahler
classes.  It was discovered by Donaldson
\cite{Do1} in the setting of moment maps and by Chen \cite{Ch2} as the gradient
flow of the $J$-functional appearing in his formula for the Mabuchi energy \cite{Ma1}.  

The J-flow is defined as follows.  Let $(M, \omega)$ be a compact K\"ahler
manifold of complex dimension $n$ and let $\chi^{\ }_0$ be another K\"ahler form on
$M$. Let $\mathcal{H}$ be the space of K\"ahler potentials
$$\mathcal{H} = \{ \phi \in C^{\infty}(M) \ | \ \chi_{\phi} = \chi^{\ }_0 +
\frac{\sqrt{-1}}{2} \partial \dbar \phi >0 \}.$$
The J-flow is the flow in $\mathcal{H}$ given by
\begin{equation} \label{eqnJflow}
\left.
\begin{array}{rcl}
{\displaystyle \ddt{\phi}} & = & \displaystyle{ c - \frac{\omega \wedge
\chi_{\phi}^{n-1}}{\chi_{\phi}^n}}
\\
\phi|_{t=0} & = & \phi_0 \in \mathcal{H}, \end{array} \right\}
\end{equation}
where $c$ is the constant given by
$$c = \frac{ [\omega] \cdot [\chi^{\ }_0]^{n-1}}{[{\chi^{\ }_0}]^n}.$$
A critical point of the J-flow gives a K\"ahler form $\chi$ satisfying
\begin{equation} \label{maineqn}
\omega \wedge \chi^{n-1} = c \chi^n.
\end{equation}

In local coordinates, this critical equation can be written
$$\chii{i}{j} \g{i}{j} = nc,$$
where $\ch{i}{j}$ and $\g{i}{j}$ are the K\"ahler metrics
corresponding to $\chi$ and $\omega$.  This is a fully
nonlinear second order elliptic equation in $\phi$. If a solution $\chi>0$ of (\ref{maineqn}) exists then it is the unique solution in its class  \cite{Ch2, Do1}.

Choosing a normal coordinate system for $g$ so that $\ch{i}{j}$ is
diagonal with entries $\lambda_1, \ldots, \lambda_n$, the critical equation
becomes
$$ \sum_{i=1}^n \frac{1}{nc \lambda_i} = 1.$$
Donaldson \cite{Do1} noted that the necessary condition
\begin{equation} \label{condition1}
\frac{1}{nc \lambda_i} < 1 \quad \textrm{for } i=1, \ldots, n,
\end{equation}
translates into the K\"ahler class condition
$$[nc {\chi^{\ }_0} - \omega ]>0.$$
He also remarked that an obvious conjecture would be that this be
sufficient for the existence of a critical metric (see
also \cite{Ch2}, Conjecture/Question 1).  Chen \cite{Ch2} confirmed this conjecture in the
case $n=2$ by reducing (\ref{maineqn}) to the complex Monge-Amp\`ere
equation which was solved by Yau \cite{Ya1}.  

For a general K\"ahler manifold, Chen \cite{Ch3} showed that solutions of
the J-flow exist for all time for any smooth initial data, and proved
convergence in the case of non-negative bisectional curvature. 
The second author showed in \cite{We1, We2} that the J-flow converges
to a critical metric under the condition 
$$nc {\chi^{\ }_0} - (n-1) \omega>0,$$
which in coordinates as above corresponds to
\begin{equation} \label{condition2}
\frac{1}{nc \lambda_i} < \frac{1}{n-1} \quad \textrm{for } i=1,\ldots,
n.
\end{equation}
This shows that the critical equation can be solved under the
class condition
$$[nc {\chi^{\ }_0} - (n-1) \omega]>0,$$
which coincides with the above conjecture for $n=2$.

In this paper we find a necessary and sufficient condition for the
convergence of the flow and existence of a critical metric.  In terms
of the $\lambda_i$, this condition can be written
\begin{equation} \label{condition3}
\sum_{\genfrac{}{}{0pt}{}{i=1}{i \neq k}}^n \frac{1}{nc \lambda_i} <1
\quad \textrm{for } k=1, \ldots, n.
\end{equation}
Clearly (\ref{condition2}) implies (\ref{condition3}) implies
(\ref{condition1}), and all three coincide for the case $n=2$.
This condition can be rewritten in terms of the positivity
of a certain $(n-1, n-1)$-form.  Precisely, we prove the following.

\addtocounter{corollary}{1}
\begin{theorem} \label{smooththeorem}  Let $M$ be a compact K\"ahler manifold of complex dimension $n$ with two K\"ahler metrics $\omega$ and ${\chi^{\ }_0}$.  The following are equivalent:
\begin{enumerate}
\item[(i)] There exists a metric $\chi'$ in $[{\chi^{\ }_0}]$
satisfying
\begin{equation} \label{condition}
(nc \chi' - (n-1)\omega) \wedge {\chi'}^{n-2}  >0.
\end{equation}
\item[(ii)]  For any initial data $\phi_0$ in $\mathcal{H}$, the $J$-flow (\ref{eqnJflow})
converges in $C^{\infty}$ to $\phi_{\infty}$ in $\mathcal{H}$ with $\chi=\chi_{\phi_{\infty}}$ satisfying the critical equation 
(\ref{maineqn}).
\item[(iii)] There exists a smooth solution $\chi$ in $[{\chi^{\ }_0}]$ to the critical equation (\ref{maineqn}).
\end{enumerate}
\end{theorem}

Recall that an $(n-1, n-1)$ form $v$ is positive if for all $(1,0)$
forms $\alpha$, 
$$ \sqrt{-1} v \wedge \alpha \wedge \overline{\alpha} >0.$$  
Notice that the condition (\ref{condition}) is satisfied if the class condition 
$$[nc {\chi^{\ }_0} - (n-1) \omega]>0$$
holds.

In \cite{Ch2}, Chen showed using results from \cite{Ch1} that finding a solution of the critical
equation implies the lower boundedness of the
Mabuchi energy \cite{Ma1} on certain K\"ahler classes for manifolds with negative
first Chern class.  We will now briefly discuss this functional and its
role in K\"ahler geometry.

The Mabuchi energy is the functional on $\mathcal{H}$ given by
$$M_{{\chi^{\ }_0}} (\phi) = - \int_0^1 \int_M \ddt{\phi_t} (R_{\chi_{\phi_t}} -
\underline{R}) \frac{\chi_{\phi_t}^n}{n!} dt,$$
where $\{\phi_t \}_{0 \le t \le 1}$ is a path in $\mathcal{H}$ between 0 and
$\phi$, $R_{\chi_{\phi_t}}$ is the scalar curvature of $\chi_{\phi_t}$ and
$\underline{R}$ is the average of the scalar curvature, given by
$$\underline{R} =  \frac{1}{\int_M
{\chi^{n}_0}}   \int_M R_{{\chi^{\ }_0}} {\chi^{n}_0}  .$$

The critical points of the Mabuchi energy are constant scalar curvature K\"ahler (cscK) metrics. 
The question of the existence of a cscK metric in a given K\"ahler
class is a difficult and interesting problem and is
 expected to be equivalent to a notion of stability in the sense
of geometric invariant theory. 
This is an idea of Yau \cite{Ya2}, who made the conjecture for K\"ahler-Einstein metrics on Fano manifolds.  In recent years there has been much progress on this problem, in particular by Tian \cite{Ti3} and Donaldson \cite{Do2}.
There is now a large body of literature pertaining to Yau's conjecture and we refer the reader to the additional references \cite{Ti2,Ti4,DiTi, Do3, Do4, Lu, PhSt1,PhSt2, PhSt3, PaTi, Ma2, RoTh}, which is far from a complete list, for details.  The Mabuchi energy and its `derivative', the Futaki invariant \cite{Fu}, are central to these ideas.  

It is known that
if there exists a K\"ahler-Einstein metric in a class $[{\chi^{\ }_0}]$ then the Mabuchi energy is bounded below on that class.  This was shown by Bando and Mabuchi \cite{BaMa} for the case $c_1(M)>0$.  The proof for the cases $c_1(M)<0$ and $c_1(M)=0$ can be found in \cite{Ti4}.    Moreover,  Chen and Tian \cite{ChTi} have recently shown that the existence of a cscK metric in any class implies the lower boundedness of Mabuchi's functional.      
In particular, a
manifold with $c_1(M)<0$ admits a K\"ahler-Einstein metric \cite{Ya1, Au1} and so the Mabuchi energy is bounded below on any K\"ahler class which is a positive multiple of $-c_1(M)$.  For manifolds with $c_1(M)=0$, there exists a Ricci-flat metric \cite{Ya1} in any class and so the Mabuchi energy is bounded below on any $[{\chi^{\ }_0}]$.

More generally, 
Tian \cite{Ti4} has conjectured that the existence of a cscK metric 
is equivalent to the `properness' of the Mabuchi energy.  This means that it is bounded below by a certain energy functional.  For the precise definition, see section 2.  This result is already known \cite{Ti3, Ti4, TiZh} when $c_1(M)$ is a multiple of $[{\chi^{\ }_0}]$.

We will deal with the case when $M$ has negative first Chern class and when $[{\chi^{\ }_0}]$ is not necessarily a multiple of $-c_1(M)$.  It was shown in \cite{We2}, using the J-flow, that if
$c_1(M)<0$ then the 
Mabuchi energy is bounded below on all K\"ahler classes $[{\chi^{\ }_0}]$
satisfying
$$ - n\frac{c_1(M) \cdot [{\chi^{\ }_0}]^{n-1}}{[{\chi^{\ }_0}]^n} [{\chi^{\ }_0}] +
(n-1)c_1(M) >0.$$
Chen \cite{Ch2} had proved the case $n=2$.  
Since this condition is easily satisfied for $[{\chi^{\ }_0}]=-c_1(M)$ it follows that this inequality defines a reasonably large open conical neighbourhood of $-c_1(M)$ in the K\"ahler cone.
Using Theorem \ref{smooththeorem} we can enlarge the set of classes for which this holds and, using a result of Tian \cite{Ti4}, show that the Mabuchi energy is not just bounded below but is in fact proper.

\addtocounter{corollary}{1}
\begin{theorem} \label{mabuchienergycorollary}
Suppose that $M$ satisfies $c_1(M)<0$. Let $V$ be the cone of all
K\"ahler classes $[{\chi^{\ }_0}]$ with the property that there exist
metrics $\omega$ in  $-\pi c_1(M)$ and 
$\chi'$ in $[{\chi^{\ }_0}]$ with
$$\left( - n\frac{\pi c_1(M) \cdot [{\chi^{\ }_0}]^{n-1}}{[{\chi^{\ }_0}]^n} \chi' - (n-1) \omega\right)
\wedge {\chi'}^{n-2}
>0.$$ Then the Mabuchi energy is proper on every class $[{\chi^{\ }_0}]$
in $V$.
\end{theorem}

In section 2 we explain the relationship between the Mabuchi energy and the J-functional which leads to this result.

Returning to the J-flow, we also consider the case when the inequality
(\ref{condition}) is not satisfied.  Using a method similar to that used by Tsuji \cite{Ts} for the K\"ahler-Ricci flow,  we can still obtain some
estimates away from a subvariety if we make more general
assumptions.  Moreover, under these assumptions, and if
(\ref{condition}) is not satisfied, we show that the flow must blow up over a subvariety in a certain sense.

\addtocounter{corollary}{1}
\begin{theorem} \label{theoremsing} Let $M$ be a compact K\"ahler manifold of complex dimension $n$ with two K\"ahler metrics $\omega$ and ${\chi^{\ }_0}$. 
Suppose that there exists an effective real divisor
$D=\sum_{\nu=1}^m a_{\nu} S_{\nu}$  on
$M$ and a metric $\chi'$ in the class $([{\chi^{\ }_0}] - (1/nc)
    [D])$
satisfying
$$(nc \chi' - (n-1)\omega) \wedge {\chi'}^{n-2} >0.$$
Let $\phi_t$ be a solution of the $J$-flow (\ref{eqnJflow}) and let $S=
\cup_{\nu=1}^m S_{\nu}$.   For any holomorphic sections $s_{\nu}$ of the line bundles associated to $S_{\nu}$,  vanishing on $S_{\nu}$, and for any hermitian metrics $h_{\nu}$ on these line bundles, there exist constants $C$, $C'$ and $A$ depending only on $D$,  $\omega$, ${\chi^{\ }_0}$, $\phi_0$,  $s_{\nu}$ and $h_{\nu}$ such that 
\begin{enumerate}
\item[(a)] 
$\displaystyle{
\phi_t(x) \ge -C + \sum_{\nu=1}^m \frac{a_{\nu}}{nc \pi} \log |s_{\nu}|^2_{h_{\nu}}(x),
  \qquad \textrm{for}
\ \, x \in M-S;}$
\item[(b)] 
$\displaystyle{\Lambda_{\omega} \chi_{\phi_t}(x) \le \frac{C'}{|s_1|_{h_1}^{2Aa_1/nc \pi} \cdots |s_m|_{h_m}^{2Aa_m/nc \pi}(x)} e^{A \phi_t(x)},
\qquad \textrm{for} \ \, x \in M-S.}$
\end{enumerate}
Let $\tilde{S}$ be the intersection of all sets $S$ corresponding to divisors in the linear system $|D|$. Then, with the same assumptions as above, if there does not exist $\chi'$ in $[{\chi^{\ }_0}]$
satisfying the condition
$$(nc \chi' - (n-1)\omega) \wedge {\chi'}^{n-2} >0,$$
then there exists a sequence of points and times $(x_i, t_i) \in M
\times [0, \infty)$ with $d(x_i, \tilde{S})
\rightarrow 0$ and $t_i \rightarrow \infty$  such that
$$(|\phi|+|\triangle_{\omega} \phi|) (x_i, t_i) \rightarrow \infty.$$
\end{theorem}

$d(x,\tilde{S})$ refers to the distance between the point $x$ and the set $\tilde{S}$ with respect to the metric $g$.

In two dimensions we will see that the conditions of this theorem are in fact always met.
Donaldson \cite{Do1} had noted
that  on K\"ahler surfaces, the condition
$$[nc {\chi^{\ }_0} - \omega]>0$$
is satisfied for \emph{all} K\"ahler classes $[{\chi^{\ }_0}]$ and $[\omega]$
if there are no curves of negative self-intersection on $M$.  He
remarked that if this inequality is violated then one might expect the
flow to blow
up over some such curves.  We confirm this conjecture in the following sense.

\addtocounter{corollary}{1}
\begin{theorem} \label{surfacetheorem}  Suppose $n
=2$ and the class 
$(nc[{\chi^{\ }_0}] - [\omega])$
is not K\"ahler.  Then there exists a positive integer $m$, irreducible curves of negative
self-intersection $E_1, \ldots, E_m$ on $M$ and positive real numbers $a_1, 
\ldots, a_m$ such that if $D= \sum_{\nu=1}^m a_{\nu} E_{\nu}$ and $S= \cup_{\nu=1}^m
E_{\nu}$ then
$$[nc {\chi^{\ }_0} - \omega] - [D]>0,$$
and we obtain the same estimates (a) and (b) of Theorem \ref{theoremsing}.  Moreover,
with $\tilde{S}$ as in that theorem, there exists a sequence of points and times $(x_i, t_i) \in M
\times [0,\infty)$ with $d(x_i, \tilde{S})
\rightarrow 0$ and $t_i \rightarrow \infty$  such that
$$(|\phi|+|\triangle_{\omega} \phi|) (x_i, t_i) \rightarrow \infty.$$
\end{theorem}

In section 2 we discuss our notation and give some preliminaries about the J-flow, the $I$ and $J$ functionals, and the Mabuchi energy.  In section 3 we prove the main estimates of Theorem \ref{theoremsing}.  Finally, in section 4 we prove the main results and make a few remarks and conjectures.

\addtocounter{section}{1}
\setcounter{equation}{0}
\bigskip
\bigskip
\noindent
{\bf 2. Preliminaries}
\bigskip

From now on, we assume that $\omega$ has been scaled so that $c=1/n.$  We
will work
in local coordinates, and write
$$\omega = \frac{\sqrt{-1}}{2} g_{i \overline{j}} dz^i \wedge
dz^{\overline{j}},
\qquad
{\chi^{\ }_0} = \frac{\sqrt{-1}}{2} \chz{i}{j} dz^i \wedge dz^{\overline{j}},$$
and
$$\chi = \frac{\sqrt{-1}}{2} \ch{i}{j} dz^i \wedge dz^{\overline{j}}
=  \frac{\sqrt{-1}}{2} (\chz{i}{j} + \dd{i} \dbr{j} \phi) dz^i \wedge
dz^{\overline{j}},$$
where $\chi= \chi_{\phi}$ (suppressing the $t$-subscript.)  The
operators
$\Lambda_{\omega}$ and
$\Lambda_{\chi}$ act on $(1,1)$ forms $\al = \frac{\sqrt{-1}}{2}
\al_{i\overline{j}} dz^i
\wedge dz^{\overline{j}}$ by
$$ \Lambda_{\omega} \al = \gi{i}{j} \al_{i \overline{j}}, \qquad
\textrm{and}
\qquad \Lambda_{\chi} \al = \chii{i}{j} \al_{i \overline{j}}.$$

The $J$-flow (\ref{eqnJflow}) can be written
\begin{equation} \label{eqnJflow2}
\left.
\begin{array}{rcl}
{\displaystyle \ddt{\phi}} & = & \displaystyle{  \frac{1}{n} (1- \Lambda_{\chi} \omega)}
\\
 \phi|_{t=0} & = & \phi_0 \in \mathcal{H}, \end{array} \right\}
\end{equation}
Differentiating with respect to $t$ gives
\begin{eqnarray} \label{eqnddt}
\ddt{}{\left(\ddt{\phi}\right)} & = & \lt \left(\ddt{\phi}\right),
\end{eqnarray}
where the operator $\lt$ acts on functions $f$ by
$$\lt f = \frac{1}{n} \h{k}{l} \dd{k} \dbr{l}
f,$$
for
$$\h{k}{l} = \chii{k}{j}\chii{i}{l}\g{i}{j}.$$

By the maximum principle for parabolic equations, (\ref{eqnddt}) implies
that
\begin{equation} \nonumber
\inf_M (\Lambda_{\chi_{\phi_0}} \omega) \le \Lambda_{\chi}\omega \le \sup_M
(\Lambda_{\chi_{\phi_0}}
\omega),
\end{equation}
which gives a lower bound for $\chi$,
\begin{equation} \label{eqnlowerbound}
\chi \ge \frac{1}{\sup_M (\Lambda_{\chi_{\phi_0}} \omega)} \,
\omega.
\end{equation}

We will now define some important functionals on the space $\mathcal{H}$.
The \emph{$J$-functional} \cite{Ch2} is defined by
$$J(\phi) = J_{\omega, {\chi^{\ }_0}} (\phi) = \int_0^1 \int_M \ddt{\phi_t} \,
\frac{\ \omega
\wedge \chi_{\phi_t}^{n-1}}{(n-1)!} \, dt,$$
where $ \{ \phi_t \}$ is a path in $\mathcal{H}$ between $0$ and $\phi$. 
The functional is independent of the choice of path.
The \emph{$I$-functional} is a well-known functional (see \cite{Ma1}) given by
\begin{equation} \nonumber
I(\phi) = I_{{\chi^{\ }_0}} (\phi) = \int_0^1 \int_M \ddt{\phi_t} \,
\frac{\chi_{\phi_t}^n}{n!} \, dt.
\end{equation}
It will also be
convenient to define a
\emph{normalized
$J$-functional} which we will denote by $\hat{J}= \hat{J}_{\omega,{\chi^{\ }_0}}$, given
by
$$\hat{J}(\phi) = J(\phi) - nc I(\phi) = \int_0^1 \int_M \ddt{\phi_t} (\omega
\wedge 
\chi_{\phi_t}^{n-1} - c \chi_{\phi_t}^n) \frac{dt}{(n-1)!}.$$ 
Note that
$\hat{J}$ has the property $\hat{J}(\phi + C) = \hat{J}(\phi)$ for constants $C$. 

The J-flow is the gradient flow of the functional $\hat{J}$ on
$\mathcal{H}$.  Alternatively, as in \cite{Ch2}, one could normalize using the $I$
functional, and consider the space
$\mathcal{H}_0 = \{ \phi \in \mathcal{H} \ | \ I(\phi) = 0 \}$. In that case,
the J-flow is the gradient flow of $J$.  We will now describe the
relationship between $\hat{J}$ and the Mabuchi energy.
Under the assumption that $\omega = - \textrm{Ric}{({\chi^{\ }_0})}>0$,
Chen's formula \cite{Ch2} for the Mabuchi energy can be written
\begin{equation} \label{mabuchienergyeqn}
M_{{\chi^{\ }_0}} (\phi) = \int_M \log \left( \frac{\chi_{\phi}^n}{{\chi^{n}_0}} \right)
\frac{\chi_{\phi}^n}{n!} + \hat{J}_{\omega, {\chi^{\ }_0}} (\phi).
\end{equation}
It is easy to see that the first term is bounded below.  In fact, we will see in section 4 that it is proper.


Let us now define what is meant by the `properness' of a functional on $\mathcal{H}$.
First, recall the definitions of the Aubin-Yau energy functionals \cite{Au2} which are also called $I$ and $J$ in the literature.  To avoid confusion we will denote them by $I^E$ and $J^E$.  They are given by
\begin{eqnarray*}
I^E_{{\chi^{\ }_0}}(\phi) & = & \frac{\sqrt{-1}}{2n! V} \sum_{i=0}^{n-1} \int_M \partial \phi \wedge \dbar \phi \wedge {\chi^{i}_0} \wedge \chi_{\phi}^{n-1-i}\\
J^E_{{\chi^{\ }_0}}(\phi) & = & \frac{\sqrt{-1}}{2n! V} \sum_{i=0}^{n-1} \frac{i+1}{n+1} \int_M \partial \phi \wedge \dbar \phi \wedge {\chi^{i}_0} \wedge \chi_{\phi}^{n-1-i},
\end{eqnarray*}
for $V=\int_M {\chi^{n}_0}/n!$, and they satisfy
$$\frac{1}{n+1} I^E_{{\chi^{\ }_0}} \le J^E_{{\chi^{\ }_0}} \le \frac{n}{n+1} I^E_{{\chi^{\ }_0}}.$$

Following Tian \cite{Ti3}, we say that a functional $T$ on $\mathcal{H}$ is \emph{proper} if there 
there exists an
increasing function
$$f: [0, \infty) \rightarrow \mathbf{R},$$
satisfying $f(x) \rightarrow \infty$ as $x \rightarrow \infty$, such that for any $\phi \in \mathcal{H}$, 
$$T(\phi) \ge f\left( J_{{\chi^{\ }_0}}^E (\phi) \right).$$

In the course of the proofs,  $C_1, C_2, \ldots$ will denote uniform constants.  Curvature expressions such as $\R{k}{l}{i}{j}$ refer to the metric $\g{i}{j}$.

\addtocounter{section}{1}
\setcounter{equation}{0}
\bigskip
\bigskip
\noindent
{\bf 3. Estimates on $\phi$ and $\Lambda_{\omega} \chi$}
\bigskip

In this section, we prove the main estimates (a) and (b) of Theorem
\ref{theoremsing}.  We assume that $c=1/n$.  As in the statement of the theorem, let $s_{\nu}$, for $\nu = 1, \ldots, m$, be holomorphic sections of the line bundles associated to the divisors $S_{\nu}$ which vanish on $S_{\nu}$ and let $h_{\nu}$ be hermitian metrics on these line bundles.  Then since $\chi' \in [{\chi^{\ }_0}] - [D]$, there exists a smooth function $\theta$ such that
$$\chi' = {\chi^{\ }_0} + \sum_{\nu =1}^m \frac{\sqrt{-1}}{2\pi} a_{\nu} \partial \dbar \log h_{\nu} + \frac{\sqrt{-1}}{2} \partial \dbar \theta.$$
Replacing $h_1$ with $h_1 e^{-\theta \pi/a_1}$, we may assume that $\theta=0$.  The change in $h_1$ will only modify the constants $C$ and $C'$ in the final estimates.

The condition for $\chi'$ gives us, for $\epsilon>0$ sufficiently small,
\begin{equation} \label{epsiloncondition}
(\chi' - (n-1)\omega) \wedge {\chi'}^{n-2} >
2\epsilon  {\chi'}^{n-1}.
\end{equation}
For $\phi$ in $\mathcal{H}$, set
$$\tphi = \phi - \sum_{\nu =1}^m \frac{a_{\nu}}{\pi}\log |s_{\nu}|^2_{h_{\nu}}.$$
Then, for any smooth $\phi$, the corresponding $\tphi$ is a smooth function on
$M-S$, and
$$\chi' = \chi_{\phi} - \frac{\sqrt{-1}}{2} \partial \dbar \tphi.$$ 

We use the maximum principle to prove the following estimate on the second
derivatives of $\phi$.

\addtocounter{theorem}{1}
\addtocounter{proposition}{1}
\begin{lemma} \label{theoremC2}
Let $\phi=\phi_t$ be a solution of the $J$-flow (\ref{eqnJflow}) on $[0,
\infty)$, with the assumptions of Theorem \ref{theoremsing}.  Then
there exist positive constants $A$ and $C_1$
depending only on $\omega$, ${\chi^{\ }_0}$, $\chi'$ and $\phi_0$ such that for any time $t\ge 0$,
$\chi =
\chi_{\phi_t}$ satisfies
\begin{equation} \label{eqnC2}
\Lambda_{\omega} \chi \le C_1 e^{A (\tphi - \inf_{(M-S) \times [0,t]} \tphi)},
\end{equation}
on $M-S$.
\end{lemma}

\begin{proof}
First note that the infimum of $\tphi$ on $M-S$ exists because $\tphi(x)$ tends to
infinity as $x$ approaches $S$.  
Choose the constant $A$ to be large enough so that
\begin{equation} \label{Ainequality}
-\frac{1}{A(\Lambda_{\omega} \chi)}(\h{k}{l} R_{k \overline{l}}^{\ \ i
\overline{j}} \ch{i}{j} - \chii{k}{l} R_{k \overline{l}}) \le \epsilon,
\end{equation}
for $\epsilon>0$ satisfying (\ref{epsiloncondition}).  Note that we can find such a constant by (\ref{eqnlowerbound}).

We will calculate the evolution of 
$$\log(\Lambda_{\omega} \chi) - A\tphi,$$
on $M-S$.  From \cite{We1}, p. 954, we have 
\begin{eqnarray*}
(\lt - \ddt{} ) \log (\Lambda_{\omega} \chi) & \ge &
\frac{1}{n(\Lambda_{\omega} \chi)}(\h{k}{l} R_{k \overline{l}}^{\ \ i
\overline{j}} \ch{i}{j} - \chii{k}{l} R_{k \overline{l}}).
\end{eqnarray*}
Calculate on $M-S$,
\begin{eqnarray*}
(\lt - \ddt{} ) \tphi & = & \frac{1}{n} ( \h{k}{l} \dd{k} \dbr{l} \tphi +
\chii{i}{j} \g{i}{j} - 1 ) \\ 
& = & \frac{1}{n} (\chii{k}{j} \chii{i}{l} \g{i}{j} \ch{k}{l} - \h{k}{l}
\chp{k}{l}  + \chii{i}{j}
\g{i}{j} - 1 )
\\ & = & \frac{1}{n} (2 \chii{i}{j} \g{i}{j} - \h{k}{l} \chp{k}{l}  - 1 ).
\end{eqnarray*}

Now consider the point $(x_0, t_0) \in (M-S) \times [0,t]$ at which the quantity
$(\log (\Lambda_{\omega}\chi) - A \tphi)$ achieves its maximum.  Note that since $\tphi(x)$ tends to infinity as $x$ approaches $S$, the maximum is achieved on this set.
  We may
assume that $t_0>0$, for if $t_0=0$, the estimate follows trivially.  At $(x_0,t_0)$, we have \begin{eqnarray*}
0 & \ge & (\lt - \ddt{} ) (\log(\Lambda_{\omega} \chi) - A \tphi) \\
& \ge & \frac{1}{n} \left( \frac{1}{(\Lambda_{\omega} \chi)}(\h{k}{l} R_{k
\overline{l}}^{\
\ i
\overline{j}} \ch{i}{j} - \chii{k}{l} R_{k \overline{l}}) 
+A\h{k}{l} \chp{k}{l} - 2A \chii{i}{j} \g{i}{j}  + A
\right).
\end{eqnarray*}
Hence at $(x_0, t_0)$,
\begin{eqnarray*}
1 + \h{k}{l} \chp{k}{l}  -
2
\chii{i}{j} \g{i}{j} 
& \le & -\frac{1}{A(\Lambda_{\omega} \chi)}(\h{k}{l}
R_{k
\overline{l}}^{\
\ i
\overline{j}} \ch{i}{j} - \chii{k}{l} R_{k \overline{l}})\\
& \le & \epsilon, 
\end{eqnarray*}
using (\ref{Ainequality}). We now compute in normal coordinates for the metric $\chi'_{i \overline{j}}$ so that the metric $\chi_{i \overline{j}}$ is diagonal with entries $\lambda_1, \ldots, \lambda_n$.  The metric $\g{i}{j}$ may not be diagonal in this basis, but we will denote its (positive) diagonal entries  $\g{i}{i}$ by $\mu_i$.  

The above
inequality becomes
$$
1 + \sum_{i=1}^n \frac{\mu_i}{\lambda_i^2}  - 2 \sum_{i=1}^n
\frac{\mu_i}{\lambda_i} \le \epsilon.$$
Fix $k$ between $1$ and $n$.
Completing the square, we obtain
\begin{equation} \label{keyinequality}
1 +  \sum_{\genfrac{}{}{0pt}{}{i=1}{i\neq k}}^n \mu_i \left( \frac{1}{\lambda_i} - 1 \right)^2 - \sum_{\genfrac{}{}{0pt}{}{i=1}{i\neq k}}^n \mu_i + \frac{\mu_k}{\lambda_k^2} - 2\frac{\mu_k}{\lambda_k} \le \epsilon.
\end{equation}
Now we will make use of the condition (\ref{epsiloncondition}).  Writing $\beta_k$ for the (1,1) form $\frac{\sqrt{-1}}{2} dz^k \wedge dz^{\overline{k}}$, we have 
$$(\chi' - (n-1)\omega) \wedge \chi'^{n-2} \wedge \beta_k > 2 \epsilon \chi'^{n-1} \wedge \beta_k.$$
Writing $\chi'$ and $\omega$ in our coordinates, this inequality becomes
$$(n-1)! \beta_1 \wedge \cdots \wedge \beta_n - (n-1)!\sum_{\genfrac{}{}{0pt}{}{i=1}{i\neq k}}^n \mu_i \beta_1 \wedge \cdots \wedge \beta_n >2\epsilon (n-1)! \beta_1 \wedge \cdots \beta_n,$$
and hence
$$1 - \sum_{\genfrac{}{}{0pt}{}{i=1}{i\neq k}}^n \mu_i >2 \epsilon.$$
Returning to (\ref{keyinequality}) this gives us
\begin{eqnarray*}
2\epsilon - 2\frac{\mu_k}{\lambda_k} & < & 1 +   \sum_{\genfrac{}{}{0pt}{}{i=1}{i\neq k}}^n \mu_i \left(
\frac{1}{\lambda_i}  - 1 \right)^2 -\sum_{\genfrac{}{}{0pt}{}{i=1}{
i \neq k}}^n \mu_i
 +\frac{\mu_k}{\lambda_k^2} - 2\frac{\mu_k}{\lambda_k}
 \le \epsilon,
\end{eqnarray*}
from which we obtain
$$\frac{\lambda_k}{\mu_k} < \frac{2}{\epsilon},$$
for $k=1, \ldots, n$.  Summing in $k$ we obtain the estimate
$$\Lambda_{\omega} \chi \le C_1 = \frac{2n}{\epsilon}$$
at
the point $(x_0, t_0).$
Then on
 on $(M-S) \times [0,t]$, we have
$$\log (\Lambda_{\omega} \chi) - A\tphi \le \log C_1 - A\inf_{(M-S) \times
[0,t]}
\tphi.$$
Exponentiating gives
$$\Lambda_{\omega} \chi \le C_1 e^{A(\tphi - \inf_{(M-S) \times [0,t]}
\tphi)},$$ completing the proof of Lemma \ref{theoremC2}.
\end{proof}

We will now prove the zero order estimate for $\phi.$

\addtocounter{theorem}{1}
\addtocounter{proposition}{1}
\begin{lemma} \label{theoremC0}
There exists a constant $C$ such that the solution $\phi=\phi_t$ of the
J-flow (\ref{eqnJflow}) satisfies
$$\displaystyle{\tphi_t(x) = \phi_t(x) - \sum_{\nu=1}^m \frac{a_{\nu}}{\pi}\log |s_{\nu}|_{h_{\nu}}^2(x)
\ge -C}, \ \textrm{for all } x \in M-S.$$
\end{lemma}
\begin{proof}
Fix a time $t$ and choose $t_0$ in the interval
$[0,t]$ so that
\begin{eqnarray*}
\inf_{M-S} \tphi_{t_0} & = & \inf_{(M-S) \times [0,t]} \tphi \\
& = &
\inf_{(M-S)\times[0,t_0]} \tphi.
\end{eqnarray*}
Define a function $\psi$ on $M-S$ by $$\psi = \tphi_{t_0} - \sup_M \phi_{t_0}.$$  
We will prove that
$\psi$ is uniformly bounded from below.   We will make use of Lemma
\ref{theoremC2}.  Set 
$$u = e^{-B\psi}, \quad  \textrm{for} \quad
B=\frac{A}{1-\delta},$$
where $A$ is the constant from Lemma \ref{theoremC2} and $\delta>0$ is a small positive constant to be
determined.  Then $u$ is a smooth non-negative function, which we will show is
uniformly bounded from above.  We have the following
lemma.

\addtocounter{proposition}{1}
\addtocounter{theorem}{1}
\begin{lemma} For any $p \ge 1$,
\begin{equation} \label{estimateu}
\int_M | \nabla u^{p/2} |^2 \frac{\omega^n}{n!} \le C_2 \, p \, \|u
\|_{C^0}^{1-\delta}
\int_M u^{p-(1-\delta)} \frac{\omega^n}{n!}.
\end{equation}
\end{lemma}

\begin{proof} The proof is a modification of the argument given in \cite{We2}.
First note that although $| \nabla u^{p/2}|^2$ may not be smooth,
it is integrable.  For any small, positive $\eta$, let $T_{\eta}$ be a tubular neighborhood of
$S$ of radius $\eta$ (with respect to the metric $g$.)  Then calculate
\begin{eqnarray*}
\int_{M-T_{\eta}} | \nabla u^{p/2}|^2 \frac{\omega^n}{n!} & = & \sqrt{-1}
\int_{M- T_{\eta}}
\partial e^{-Bp\psi/2}
\wedge \dbar e^{-Bp\psi/2} \wedge \frac{\omega^{n-1}}{(n-1)!} \\
& = & \frac{B^2 p^2}{4} \sqrt{-1} \int_{M-T_{\eta}} e^{-Bp\psi} \partial \psi
\wedge
\dbar
\psi \wedge \frac{\omega^{n-1}}{(n-1)!} \\
& = & - \frac{Bp}{4} \sqrt{-1} \int_{M-T_{\eta}} \partial (e^{-Bp\psi}) \wedge
\dbar
\psi
\wedge \frac{\omega^{n-1}}{(n-1)!} \\
& = & \frac{Bp}{2} \int_{M-T_{\eta}} e^{-Bp\psi} \frac{\sqrt{-1}}{2} \partial
\dbar
\psi \wedge \frac{\omega^{n-1}}{(n-1)!} + E_{\eta},
\end{eqnarray*}
where $E_{\eta}$ is the boundary term,
$$E_{\eta} = - \frac{Bp}{4} \sqrt{-1} \int_{\partial(M-T_{\eta})} e^{-Bp\psi} \,
\dbar
\psi
\wedge \frac{\omega^{n-1}}{(n-1)!},$$
obtained by integrating by parts. Note that
$E_{\eta} \rightarrow 0$ as $\eta \rightarrow 0$.  Now from the definition of
$\psi$,
\begin{eqnarray*}
\int_{M-T_{\eta}} | \nabla u^{p/2}|^2 \frac{\omega^n}{n!} & = &
\frac{Bp}{2}
\int_{M-T_{\eta}} e^{-Bp\psi} (\chi_{\phi_{t_0}} - \chi')
\wedge
\frac{\omega^{n-1}}{(n-1)!} + E_{\eta}\\ & \le & \frac{Bp}{2} \int_{M-T_{\eta}}
e^{-Bp
\psi} (\Lambda_{\omega} \chi_{\phi_{t_0}}) \frac{\omega^n}{n!} + E_{\eta} \\
& \le & \frac{C_1 Bp}{2} \int_{M-T_{\eta}} e^{-Bp\psi} e^{A(\psi - \inf_{M-S} \psi)}
\frac{\omega^n}{n!} + E_{\eta} \\ 
& = & \frac{C_1 Bp}{2} e^{-A\inf_{M-S}{\psi}}
\int_{M-T_{\eta}} e^{-(p-(1-\delta))B\psi}
\frac{\omega^n}{n!} +E_{\eta} \\
& \le & \frac{C_1 Bp}{2} \| u \|_{C^0}^{1-\delta} \int_M u^{p-(1-\delta)}
\frac{\omega^n}{n!} + E_{\eta},
\end{eqnarray*}
where, in the third line, we have used the estimate
$$\Lambda_{\omega} \chi \le C_1 e^{A( \tphi_{t_0} - \inf_{(M-S) \times [0,t_0]} \tphi)} =
C_1 e^{A(\psi - \inf_{(M-S)} \psi)},$$
of Theorem \ref{theoremC2}.  Letting $\eta \rightarrow 0$ completes the proof.
\end{proof}

We now use the following lemma from \cite{We2}, which we quote without proof.

\addtocounter{proposition}{1}
\addtocounter{theorem}{1}
\begin{lemma}
If $u \ge 0$ satisfies the estimate (\ref{estimateu}) for all $p \ge 1$, then for
some constant $C_3$ independent of $u$,
$$ \| u \|_{C^0} \le C_3\left( \int_M u^{\delta}
\frac{\omega^n}{n!} \right)^{1/\delta}.$$
\end{lemma}

To obtain the upper bound for $u$, observe that
\begin{eqnarray*}
\int_M u^{\delta}\frac{\omega^n}{n!} & = & \int_M |s_1|_{h_1}^{2B\delta  a_1/ \pi} \cdots |s_m|_{h_m}^{2B\delta a_m /\pi}
e^{-B\delta(\phi_{t_0} -
\sup_M
\phi_{t_0})}\frac{\omega^n}{n!} \\
& \le & C_4 \int_M e^{-B\delta (\phi_{t_0} - \sup_M
\phi_{t_0})}\frac{{\chi^{n}_0}}{n!}.
\end{eqnarray*}
Choosing $\delta$ small enough we can bound the right hand side.  This is due to the following
proposition of Tian \cite{Ti1}, based on a result of H\"ormander \cite{Ho} (it was used by Tian to define the $\alpha$-invariant - see \cite{TiYa, So} for more on this important invariant.)

\addtocounter{lemma}{1}
\addtocounter{theorem}{1}
\begin{proposition} \label{propositionalpha}
There exists $\alpha>0$ and $C_5$ depending only on $(M, {\chi^{\ }_0})$ such that
\begin{equation} \label{eqnalpha}
\int_M e^{-\alpha \phi} \frac{{\chi^{n}_0}}{n!} \le C_5,
\end{equation}
for all $\phi \in C^2(M)$ satisfying
$$ {\chi^{\ }_0} + \frac{\sqrt{-1}}{2} \partial \dbar \phi >0, \qquad \sup_M \phi = 0.$$
\end{proposition}

Hence we obtain an upper bound for $u$, and so a lower bound for $\psi$.  To
get the lower bound for $\tphi$, we require the following lemma from
\cite{We2}.

\addtocounter{proposition}{1}
\addtocounter{theorem}{1}
\begin{lemma}  \label{lemmasupinf} Let $\phi_t$ be a solution of the J-flow.  There exist positive constants $C_6$
  and $C_7$ depending only on the initial data such that
$$-C_6 \le \sup_M \phi_t \le C_6 - C_7 \inf_M \phi_t.$$
\end{lemma}

Using the first inequality we obtain at time $t$
\begin{eqnarray*}
\inf_{M-S} \tphi_t & \ge & \inf_{M-S} \tphi_{t_0} \\
& = & \inf_{M-S} \psi + \sup_M \phi_{t_0} \\
& \ge & -C_8,
\end{eqnarray*}
completing the proof of the lower bound of $\tphi$ and giving (a).  Estimate (b) follows immediately from (a) and Lemma \ref{theoremC2}.
\end{proof}

\addtocounter{section}{1}
\setcounter{lemma}{0}
\setcounter{proposition}{0}
\setcounter{theorem}{0}
\setcounter{equation}{0}
\bigskip
\bigskip
\noindent
{\bf 4. Proofs of the main results}
\bigskip

In section 3 we proved the estimates (a) and (b) of
Theorem \ref{theoremsing}.  We will use these to complete the proofs
of the remaining theorems.

\bigskip
\noindent
{\it \bf  Proof of Theorem \ref{smooththeorem}.}
(ii) implies (iii) is trivial.  
To see that (iii) implies (i), suppose that $\chi$ is a solution of the critical equation and choose a normal coordinate system for $\g{i}{j}$ in which $\ch{i}{j}$ is diagonal with entries $\lambda_1, \ldots, \lambda_n$.  
Assuming that $c=1/n$, we have, as explained in section 1, the inequality
$$\sum_{\genfrac{}{}{0pt}{}{i=1}{i \neq k}}^n \frac{1}{\lambda_i} <1
\quad \textrm{for } k=1, \ldots, n.$$
We will now calculate the form $(\chi - (n-1)\omega)\wedge \chi^{n-2}$.
As in section 3, write $\beta_k$ for the (1,1)-form $\frac{\sqrt{-1}}{2} dz^k \wedge d\overline{z}^k$.  Then
\begin{eqnarray*}
\chi^{n-1} & = & \left( \lambda_1 \beta_1 + \cdots  + \lambda_n \beta_n \right)^{n-1} \\
& = & (n-1)! \sum_{k=1}^n \lambda_1 \cdots \hat{\lambda_k} \cdots \lambda_n \, \beta_1 \wedge \cdots \wedge \hat{\beta_k} \wedge \cdots \wedge \beta_n, \\
\end{eqnarray*}
where `$\ \hat{} \ $' indicates that the symbol should be omitted.  Similarly,
\begin{eqnarray*}
\omega \wedge \chi^{n-2} & = & \left( \beta_1 + \cdots + \beta_n \right) \wedge \left(  \lambda_1 \beta_1 + \cdots + \lambda_n \beta_n  \right)^{n-2} \\
& = & (n-2)! \sum_{k=1}^n \sum_{\genfrac{}{}{0pt}{}{i=1}{i \neq k}}^n   \lambda_1 \cdots \hat{\lambda_i} \cdots \hat{\lambda_k} \cdots \lambda_n \beta_1 \wedge \cdots \wedge \hat{\beta_k} \wedge \cdots \wedge \beta_n.
\end{eqnarray*}
Then the condition 
$$(\chi - (n-1) \omega) \wedge \chi^{n-2} >0,$$
is equivalent to
$$ \lambda_1 \cdots \hat{\lambda_k} \cdots \lambda_n - \sum_{\genfrac{}{}{0pt}{}{i=1}{i \neq k}}^n \lambda_1 \cdots \hat{\lambda_i} \cdots \hat{\lambda_k} \cdots  \lambda_n >0, \qquad \textrm{for } k=1, \ldots, n,$$
which, since the $\lambda_i$ are positive, is precisely our inequality
$$\sum_{\genfrac{}{}{0pt}{}{i=1}{i \neq k}}^n \frac{1}{\lambda_i} <1
\quad \textrm{for } k=1, \ldots, n.$$

To show that (i) implies (ii) we will use the estimates (a) and (b) of Theorem \ref{theoremsing}.  Let $D$ be the zero divisor and let $s$ be equal to
the constant section 1.  Then estimate (a) gives us 
a uniform lower bound for $\phi$ along
the flow.  We can apply Lemma \ref{lemmasupinf} to obtain a uniform upper
bound for $\phi$.  Estimate (b) gives us a uniform bound on
$\Lambda_{\omega} \chi_{\phi}$ and hence on the second derivatives $\dd{i}
\dbr{j} \phi$.  Bounds on all the derivatives of $\phi$ and the
convergence to a critical metric then follow by the same arguments as in \cite{Ch3,We1,We2}.  This completes the proof.

\bigskip
\noindent
{\it \bf Proof of Theorem \ref{mabuchienergycorollary}.}
Arguing as in \cite{Ch2}, we can apply Yau's Theorem \cite{Ya1} and
assume that $-\omega = \textrm{Ric}({\chi^{\ }_0})$.   Since the conditions of Theorem
\ref{smooththeorem} are satisfied, the $J$-flow converges to a smooth critical
metric.  Since the $J$-flow is the gradient flow for $\hat{J}_{\omega,
  {\chi^{\ }_0}}$ and the critical metrics are unique  
it follows that $\hat{J}_{\omega, {\chi^{\ }_0}}$ is bounded below (alternatively, apply Proposition 3 of \cite{Ch2}.)
Then from
(\ref{mabuchienergyeqn}), we have
$$M_{{\chi^{\ }_0}}(\phi) \ge \int_M \log \left( \frac{\chi_{\phi}^n}{{\chi^{n}_0}} \right)
\frac{\chi_{\phi}^n}{n!} - C_1.$$

The theorem then follows immediately from a lemma due to Tian \cite{Ti4} and the properties of $I^E_{{\chi^{\ }_0}}$ and $J^E_{{\chi^{\ }_0}}$.

\addtocounter{proposition}{1}
\addtocounter{remark}{1}
\begin{lemma}
There exist positive constants $\alpha$ and $C_3$ such that
$$\frac{1}{V} \int_M \log \left( \frac{\chi_{\phi}^n}{{\chi^{n}_0}} \right)
\frac{\chi_{\phi}^n}{n!} \ge \alpha I^E_{{\chi^{\ }_0}}(\phi) - C_3.$$
\end{lemma}
\begin{proof}
The result can be found in \cite{Ti4}, p.95 (it is assumed there that $c_1(M)$ is a multiple of $[{\chi^{\ }_0}]$, but this plays no part in the proof.)  For the reader's convenience, we will give the proof here.
From Proposition \ref{propositionalpha} there exist positive constants $C_2$ and $\alpha$ such that
\begin{eqnarray*}
\frac{1}{V} \int_M e^{- \log \left(  \frac{\chi_{\phi}^n}{{\chi^{n}_0}} \right) - \alpha (\phi - \sup_M \phi)} \frac{\chi_{\phi}^n}{n!} & = & \frac{1}{V} \int_M e^{- \alpha (\phi - \sup_M \phi )} \frac{{\chi^{n}_0}}{n!} \\
& \le & C_2.
\end{eqnarray*}
By the convexity of the exponential function, we have
$$\frac{1}{V} \int_M \log \left( \frac{\chi_{\phi}^n}{{\chi^{n}_0}} \right)
\frac{\chi_{\phi}^n}{n!} \ge - \frac{\alpha}{V} \int_M (\phi - \sup_M \phi) \frac{\chi_{\phi}^n}{n!}   - \log C_2.$$
Writing $I^E_{{\chi^{\ }_0}}$ in the form 
$$I^E_{{\chi^{\ }_0}}(\phi)  = \frac{1}{n! V} \int_M \phi ({\chi^{n}_0} - \chi_{\phi}^n),$$
completes the proof of the lemma and hence the theorem.
\end{proof}


\addtocounter{proposition}{1}
\begin{remark} From Tian's conjecture \cite{Ti4} on the equivalence of the existence of a cscK metric and the properness of the Mabuchi energy, we would expect that there exists a cscK metric in each class $[{\chi^{\ }_0}]$ in the set $V$ described in Theorem \ref{mabuchienergycorollary} (c.f. \cite{We3}, Conjecture 5.2.1.)  Moreover, on the boundary of $V$, where $\chi'$ and $\omega$ satisfy
$$\left( - n\frac{\pi c_1(M) \cdot [{\chi^{\ }_0}]^{n-1}}{[{\chi^{\ }_0}]^n} \chi' - (n-1) \omega\right)
\wedge {\chi'}^{n-2}
\ge 0,$$
we would expect the Mabuchi energy to be bounded below.  Also, one would expect the classes in $V$ to be K-stable and those on the boundary to be K-semistable \cite{Ti3}. 
\end{remark}

\bigskip
\noindent
{\it \bf Proof of Theorem \ref{theoremsing}.}
It is only left to prove the last statement of this theorem.  Suppose
for a contradiction that there exists some $\eta>0$ such that
$$(|\phi| + |\triangle_{\omega} \phi|)(x) \le C_4, \qquad \textrm{for } x \in T_{\eta},$$
where $T_{\eta}$ is a tubular neighborhood of $\tilde{S}$ of radius $\eta$.  We will
show this implies that the J-flow
converges to a critical metric, contradicting our assumption that there
does not exist $\chi'$ in $[{\chi^{\ }_0}]$ satisfying (\ref{condition}).

First, it is not difficult to see that there exist a finite number of divisors $D^{(1)}, \ldots, D^{(k)}$ in the linear system $|D|$ such that any $x$ in $M-\overline{T}_{\eta}$ is at least a distance $\eta/2$ from one of the $D^{(i)}$.

We claim that there exists $C_5$ such that for $i=1, \ldots, k$,
\begin{eqnarray*}
\inf_{M-S^{(i)}} \tphi^{(i)} 
& \le & C_5,
\end{eqnarray*}
where we are using the obvious notation. 
To see this, pick for each $i$ a point $x^{(i)}$ in $T_{\eta}$ which does not lie on  $S^{(i)}$.  Then
\begin{eqnarray*}
\inf_{M-S^{(i)}} \tphi^{(i)} & \le & \phi(x^{(i)}) - \sum_{\nu=1}^{m^{(i)}} \frac{a^{(i)}_{\nu}}{\pi} \log|s^{(i)}_{\nu} |^2_{h^{(i)}_{\nu}}(x^{(i)}),
\end{eqnarray*}
which is uniformly bounded.
Then from estimate (b) of Theorem \ref{theoremsing}, we have for $i=1, \ldots, k$,
$$ \Lambda_{\omega} \chi \le C_6 e^{A(\tphi^{(i)} - \inf_{M-S^{(i)}} \tphi^{(i)})}.$$
Note that this is a stronger estimate than (\ref{eqnC2}), and
enables us to do the following.  Define $\psi^{(i)} = \tphi^{(i)} - \sup_M \phi$
as before.  Then we have \emph{at any time $t$},
$$\Lambda_{\omega} \chi \le C_6 e^{A(\psi^{(i)} - \inf_{M-S^{(i)}} \psi^{(i)})}.$$
By the same argument as in section 3, we obtain a uniform constant
$C_7$ such that
$$\psi^{(i)} \ge - C_7.$$
Hence on $M-S^{(i)}$,
\begin{eqnarray*}
\Lambda_{\omega} \chi & \le & C_6 e^{A(\psi^{(i)} - \inf_{M-S^{(i)}} \psi^{(i)})} \\
& \le & 
C_8 e^{A \psi^{(i)}} \\
& = & C_8 e^{A(-\sum_{\nu=1}^{m^{(i)}} \frac{a^{(i)}_{\nu}}{\pi} \log|s^{(i)}_{\nu} |^2_{h^{(i)}_{\nu}} + \phi - \sup_M \phi)} \\ 
& \le & \frac{C_8}{|s^{(i)}_1|^{2Aa^{(i)}_1/\pi}_{h^{(i)}_1} \cdots |s^{(i)}_{m^{(i)}}|^{2Aa^{(i)}_{m^{(i)}} /\pi}_{h^{(i)}_{m^{(i)}}}  }.
\end{eqnarray*}
Then we see that $\Lambda_{\omega} \chi$ is uniformly bounded on $M - \overline{T}_{\eta}$
and therefore on $M$.  Arguing now as in the smooth case, we obtain
a uniform bound on $\phi$ and so the J-flow converges to a critical
metric, giving us the contradiction.

\addtocounter{proposition}{1}
\begin{remark}
It is not difficult to see by the above arguments that the set $\tilde{S}$ of Theorem \ref{theoremsing} is non-empty unless
 (\ref{condition}) is satisfied for some metric $\chi'$ in $[\chi_0]$. 
\end{remark}

\addtocounter{proposition}{1}
\begin{remark}
It would be interesting to know whether one can improve on these estimates (for example, by showing that $\sup \phi$ is bounded, at least away from the singular set.)
We conjecture that in the boundary case, when there exists a metric $\chi'$ in $[{\chi^{\ }_0}]$ satisfying
$$(nc \chi' - (n-1)\omega) \wedge {\chi'}^{n-2}  \ge 0,$$
the J-flow converges to a critical metric on compact subsets outside the singular set $S$. 
\end{remark}

\bigskip
\noindent
{\it \bf Proof of Theorem \ref{surfacetheorem}.}
This theorem follows almost immediately from the following
proposition.

\addtocounter{remark}{1}
\begin{proposition}
Let $M$ be a K\"ahler surface with a K\"ahler class $\beta$ in
$H^{1,1}(M, \mathbf{R})$.   If $\alpha$ in $H^{1,1}(M, \mathbf{R})$
satisfies
$$\alpha^2 >0 \quad \textrm{and} \quad \alpha \cdot \beta >0$$
then either $\alpha$ is a K\"ahler class or
there exists a positive integer $m$, curves of negative self
intersection $E_1, \ldots, E_m$ and positive real numbers $a_1,
\ldots, a_m$ such that
$$\alpha -\sum_{\nu=1}^m a_{\nu} [E_{\nu}] $$
is a K\"ahler class.
\end{proposition}

\begin{proof}
This result is essentially contained in \cite{La} (see also \cite{Bu}) and so we 
will just give an outline of the proof here.  By Lemma 5.2 and Theorem 5.1 of \cite{La}, the conditions 
$\alpha^2>0$ and $\alpha \cdot \beta>0$ imply the existence of a K\"ahler current $\tau$ such that
$\alpha = [\tau]$.  That is, $\tau$ is a closed (1,1) current satisfying $\tau \ge \psi$ for some strictly positive smooth (1,1) form $\psi$.  By Siu's decomposition \cite{Si} and a result of Demailly \cite{De}, there exist constants $c_{\nu} \ge 0$ and irreducible curves $D_{\nu}$ such that
$$\tau = \overline{\tau} + \sum_{\nu=1}^{\infty} c_{\nu} D_{\nu},$$
where $\overline{\tau}$ is a K\"ahler current which is smooth away from a finite number of points.   By a smoothing argument, the class $[\overline{\tau}]$ is K\"ahler.
Write for each $\nu$, $$c_{\nu} D_{\nu} = a_{\nu} E_{\nu} + b_{\nu} C_{\nu},$$
where the $E_{\nu}$ and $C_{\nu}$ are irreducible curves satisfying $E_{\nu}^2 \le 0$ and $C_{\nu}^2 \ge 0$ and $a_{\nu}$ and $b_{\nu}$ are nonnegative constants.  We have
$$[\tau - \sum_{\nu=1}^m a_{\nu} E_{\nu}] = [\overline{\tau} + \sum_{\nu=1}^{\infty} b_{\nu} C_{\nu}] + [\epsilon_m],$$
where $\epsilon_m$ is the current
$$\epsilon_m = \sum_{\nu=m+1}^{\infty} a_{\nu} E_{\nu},$$
which tends to zero (in the weak topology of currents) as $m$ tends to infinity.  Now the $C_{\nu}$ are nef  and the K\"ahler cone is stable under the addition of nef classes \cite{La,Bu} so the class $[\overline{\tau} + \sum_{\nu=1}^{\infty} b_{\nu} C_{\nu}]$ is K\"ahler.   Since the K\"ahler cone is  open, there exists $m$ large enough such that
$$[\tau - \sum_{\nu=1}^m a_{\nu} E_{\nu}] = \alpha - \sum_{\nu=1}^m a_{\nu} [E_{\nu}]$$
is K\"ahler.
\end{proof}

We can apply this proposition in our case to $\alpha = [nc {\chi^{\ }_0} - \omega]$, since
$$[nc {\chi^{\ }_0} - \omega]^2 = [\omega]^2 >0,$$
and 
$$[nc {\chi^{\ }_0} - \omega] \cdot [{\chi^{\ }_0}] = [\omega] \cdot [{\chi^{\ }_0}] >0.$$
The rest of the theorem follows immediately from Theorem \ref{theoremsing}.

\begin{remark}
Can Theorem \ref{surfacetheorem} be generalized to higher dimensions?
\end{remark}

\begin{remark} If a surface $M$ with $c_1(M)<0$ has no curves of negative self-intersection then from Theorem \ref{mabuchienergycorollary}  we see that the Mabuchi energy is proper for any class.  It would be interesting to see examples of such manifolds where cscK metrics can be constructed in those classes away from the canonical class (c.f. 
 \cite{Fi}.)  If there do exist curves of negative self intersection, then one might guess that they form obstructions to the stability of $(M, [{\chi^{\ }_0}])$ in some sense (this could be related to the slope stability of \cite{Ro}, \cite{RoTh}.)

\end{remark}

\bigskip
\noindent
{\bf Acknowledgements}
The authors would like to thank:  Professor D.H. Phong, their thesis advisor, for his continued support, encouragement and advice; Professor J. Sturm, for some useful discussions and in particular for suggesting a condition of the form of (\ref{condition}); Professor G. Tian, for his advice and helpful suggestions; and Professor S.-T. Yau, for his encouragement and for some enlightening discussions.
Part of this work was carried out while the authors were supported by graduate fellowships at Columbia University.

\end{document}